\documentclass[12pt,twoside]{article}

\headheight=\baselineskip
\usepackage{amsthm,amsmath,amsfonts}
\usepackage[logical-quotes]{amsrefs}

\DeclareMathOperator{\SO}{SO}
\DeclareMathOperator{\Cl}{Cl}
\DeclareMathOperator{\symplectic}{Sp}
  \def\Sp{\symplectic}
\DeclareMathOperator{\GL}{GL} 
\DeclareMathOperator{\FullOrthog}{O} 
  \renewcommand{\O}{\FullOrthog}

\newcommand{\Z}{\mathbb{Z}}
\renewcommand{\k}{k}
\renewcommand{\phi}{\varphi}

\newcommand{\rtsys}{\Phi}
\newcommand{\rtbase}{\Delta}
\newcommand{\size}[1]{|#1|}
\newcommand{\double}[3]{{#1}{\backslash}{#2}{/}{#3}}
\newcommand{\sdouble}[3]{\size{\double{#1}{#2}{#3}}}
\newcommand{\st}{{\ \mid \ }}

\newcommand{\onto}{\to\kern-.7em\to}
\newcommand{\genby}[1]{\langle #1 \rangle}
\renewcommand{\tilde}{\widetilde}
\renewcommand{\hat}{\widehat}
\newcommand{\perpplus}{\perp}

\newcommand{\inv}{{}^{-1}}

\newcommand{\compress} {\setlength{\partopsep}{0in}
  \setlength{\itemsep}{0in} 
  \setlength{\parsep}{0in}\setlength{\parskip}{0in}}
  
\newcommand{\romanate}
{\renewcommand{\labelenumi}{{\rm(\roman{enumi})}}
 \renewcommand{\labelenumii}{{\rm(\alph{enumii})}}}

\newcommand{\alpherate}
{\renewcommand{\labelenumi}{{\rm(\alph{enumi})}}
 \renewcommand{\labelenumii}{{\rm(\alph{enumii})}}}
 
\newenvironment{theorem_count}
  {\begin{enumerate}\compress\romanate}{\end{enumerate}} 

\newcommand{\Theorempart}[2]{Theorem~{\rm\ref{#1}~\!#2}}
\newcommand{\Lemmapart}[2]{Lemma~{\rm\ref{#1}~\!#2}}
\newcommand{\Corpart}[2]{Corollary~{\rm\ref{#1}~\!#2}}

\newenvironment{map_machine}[5]
   {
    
    \begin{array}[#1]{rcl}
   #2 & {-}\!\!{-}\!\!{\longrightarrow} & #3\\ 
   #4 & \mapstochar\kern-.25em{-}\!\!{\longrightarrow} & #5
   }{\end{array}}

\newtheorem{lemma}{Lemma}[section]
\newtheorem{theorem}{Theorem}
\newtheorem{cor}[lemma]{Corollary}
\newtheorem{lemma_case}[lemma]{Lemma}

\theoremstyle{definition}
\newtheorem{remark}[lemma]{Remark}

\newtheorem{InnerExample}[lemma]{Example}
\newenvironment{example}[1]{\begin{InnerExample} #1.}{\end{InnerExample}}

\newtheorem{InnerRemarks}[lemma]{Remarks}
\newcounter{RemarkCounter}
\newcommand{\rem}{\stepcounter{RemarkCounter}{\rm 
    (\alph{RemarkCounter})}\ }
\newenvironment{remarks}{\setcounter{RemarkCounter}
             {0}\begin{InnerRemarks}}{\end{InnerRemarks}}

\newcommand{\parlabel}{}
\newtheorem*{innerpara}{\parlabel}

\title{A Classification of Certain Finite Double Coset Collections in
  the Classical Groups}

\author{W. Ethan
  Duckworth\footnote{\emph{AMS subject
      codes:} 14L, 20G15. \emph{Keywords:} Classical groups, double
    cosets, finite orbits}}

\date{May 23, 2003}

\begin{document}
\maketitle
\begin{abstract}
  Let $G$ be a classical algebraic group, $X$ a maximal rank reductive
  subgroup and $P$ a parabolic subgroup.  This paper classifies when
  $\double XGP$ is finite.  Finiteness is proven using geometric
  arguments about the action of $X$ on subspaces of the natural module
  for $G$.  Infiniteness is proven using a dimension criterion which
  involves root systems.
\end{abstract}

\section{Statement of results}
\label{section_statetement_of_results}
Let $G$ be a classical algebraic group defined over an algebraically
closed field, let $X$ be a maximal rank reductive subgroup, and let
$P$ be a parabolic subgroup.  The property of finiteness for $\double
XGP$ is preserved under taking isogenies, quotients by the center of
$G$, connected components and conjugates (see
Lemma~\ref{lemma_prelimnary_reductions} for a precise statement).
Thus, if desired, we can specify only the Lie type of $G$.  Similarly,
we can specify only the conjugacy class of $X$ and $P$; thus we
usually give the Lie type of $X$ and describe $P$ by crossing off
nodes from the Dynkin diagram for $G$.  For the purpose of classifying
finiteness, it suffices to consider only those $X$ which are defined
over $\Z$.

A subgroup $X$ is spherical if $\double XGB$ is finite for some (or,
equivalently, for each) Borel subgroup $B$.  For each classical group
we list in Table~\ref{table_spherical_subgroups} those maximal rank
reductive spherical subgroups which are defined over $\Z$.
\begin{table*}
\caption{Maximal rank reductive spherical subgroups defined over $\Z$}%
\label{table_spherical_subgroups}%
$$\begin{array}{rclrcl}
\hline
X &\le&   G  & \qquad X & \le&  G\qquad \\[.5\jot]\hline\hline
\rule{3ex}{0in} 
A_nA_mT_1  & \le &  A_{n+m+1} &  B_nD_m     & \le &  B_{n+m} \\
C_nC_m     & \le &  C_{n+m}   &  A_{n-1}T_1 & \le &  B_n     \\  
C_{n-1}T_1 & \le &  C_n       & D_nD_m     & \le &  D_{n+m}  \\
A_{n-1}T_1 & \le &  C_n       & A_{n-1}T_1 & \le &  D_n      \\
\hline
\end{array}$$
\end{table*}
We first describe the notation which is used for the list, and for the
rest of the paper, and then describe how the list is obtained.  We
write $X=A_{n} A_{m}T_1$ if $X$ is a group of Lie type $A_{n} + A_{m}$
which has a $1$-dimensional central torus, and we use similar notation
for other subgroups.  If $G$ equals $D_n$ we adopt a notational
convention to distinguish between certain subgroups of the same Lie
type which are not conjugate.  In $G = \SO(V)$ any factor denoted by
$D_{n_1}$ (or $\SO_{2n_1}$) acts as $\SO(V_1)$ for some decomposition
$V= V_1\perp V_2$ and any two factors denoted by $A_{n_1}T_1$ (or
$\GL_{n_1+1}$) act as $\GL_{n_1+1}$ on a pair of totally singular
subspaces $E$ and $F$ such that $V = (E\oplus F) \perp V_2$ (in
particular $\dim E = \dim F = n_1+1$ and $E$ and $F$ are in duality).
We allow the notation $A_0$, $B_0$ and $C_0$ to denote trivial groups
and $D_1$ to denote a group which is a $1$-dimensional torus.  We now
describe how Table~\ref{table_spherical_subgroups} is obtained.
Kr\"amer \cite{kramer} classified the reductive spherical subgroups in
characteristic $0$.  The subgroups on Kr\"amer's list were shown to be
spherical in all characteristics by Brundan \cite{brundan} and Lawther
\cite{lawther}.  Duckworth \cite{duckworth} showed that this list is
complete for maximal rank subgroups.

In Theorem~\ref{classification} we use the notational conventions just
described, as well as the following.  We write $X=L_i$ if
$X$ is conjugate to a Levi subgroup obtained by crossing off node $i$
from the Dynkin diagram of $G$ (we number the nodes of the Dynkin
diagram as in \cite{bourbaki}).  The meaning of $X=L_{i_1,\,i_2}$ and
$P=P_i$ is similar.

%
%
%
%
%
%

\begin{theorem}
\label{classification}
Let $G$ be a simple algebraic group of type $A_n$, $B_n$, $C_n$ or
$D_n$, let $X$ be a maximal rank reductive subgroup defined over $\Z$
and $P\ne G$ a parabolic subgroup.  Then $\double XGP$ is finite if
and only if $X$ is spherical or one of the following holds:
\vspace{-\topsep}
\begin{enumerate}
\compress\romanate

\item $G=A_n$,
\vspace{-.5\topsep}
\begin{enumerate}
\compress\alpherate 

\item $P\in \{P_1,\ P_n\}$,

\item $X = A_{n_1} A_{n_2} A_{n_3} T_2$ and $P = P_i$ for some $i$;
\end{enumerate}
\vspace{-.5\topsep}
\item $G=B_n$, $X = A_{n_1}B_{n_2}T_1$ and $P \in \{P_1,\ P_n\}$;

\item $G=C_n$,
\vspace{-.5\topsep}
\begin{enumerate}
\compress\alpherate 

\item $X\in \{C_{n_1}\cdots C_{n_r}$,\ \ $A_{n_1} C_{n_2}\cdots
C_{n_r}T_1\}$ and $P=P_1$,

\item $X\in \{C_{n_1}C_{n_2}C_{n_3}$,\ \  $C_{n_1}C_{n_2}T_1$,\ \ 
$A_{n_1}C_{n_2}T_1\}$ and $P=P_n$;

\end{enumerate}
\vspace{-.5\topsep}
\item $G=D_n$,
\vspace{-.5\topsep}
\begin{enumerate}
\compress\alpherate

\item $G=D_4$, $(X,P) \in \{(L_{2,3}, P_4),\ (L_{2,4},\ 
  P_3)\}$,


  
  
\item $X \in \{A_{n_1}D_{n_2}T_1$,\ \ $A_{n_1}A_{n_2}T_2\}$
and $P = P_1$,

\item $X\in \{D_{n_1}D_{n_2}D_{n_3}$,\ \  $A_{n_1}D_{n_2}T_1\}$ and
$P\in\{P_{n-1}$,\ \  $P_n\}$.
\end{enumerate}
\end{enumerate}
\end{theorem}

\section{History and preliminaries}
The problem of identifying finite double coset collections has been
studied by a variety of authors, usually in the context of some other
problem.  The spherical subgroups discussed above provide one example
of this study.  We briefly describe another example here, and refer
the reader to \cite{duckworth} and \cite{seitz} for fuller
discussions.

\begin{example}{We describe irreducible finite orbit modules}
  Let $V$ be a finite dimensional vector space and let $X$ be a
  closed, connected subgroup of $G = \GL(V)$ such that $X$ has
  finitely many orbits on $V$ and $V$ is irreducible under $X$.  Such
  modules were classified by Kac \cite{kac} in characteristic $0$ and
  by Guralnick, Liebeck, Macpherson and Seitz \cite{GLMS} in positive
  characteristic.  This work was representation theoretic, however it
  is related to the double coset question studied in the present
  paper.  Note that $X$ is a reductive group and that $\double
  X{\GL(V)}{P_1}$ is finite where $P_1$ is the stabilizer of a
  $1$-space.  The authors of \cite{GLMS} mentioned this connection
  with double cosets, they classified finiteness for $\double
  X{\GL(V)}{P_i}$ where $1\le i \le \dim V$ and they established
  \Theorempart{classification}{(i)} using rather different arguments
  from those that appear in the present paper.
\end{example}

The following lemma provides basic reductions in our double coset
question.

\begin{lemma}
\label{lemma_prelimnary_reductions}
Let $G$ be a group and let $X$ and $P$ be subgroups.  Let $Z$ be the
center of $G$, suppose that $Z\le P$ and let $\overline X$, $\overline G$
and $\overline P$ be the images of $X$, $G$ and $P$, respectively,
under the map $G\to G/Z$.  Let $K$ be a finite normal subgroup of $G$
and let $\hat X$, $\hat G$ and $\hat P$ be the images of $X$, $G$ and
$P$, respectively, under the map $G\to G/K$.  Let $g,h \in G$.  The
following are equivalent:
\begin{theorem_count}
\item $\sdouble XGP<\infty$,
\item $\sdouble {\hat X}{\hat G}{\hat P}<\infty$,
\item $\sdouble {\overline X}{\overline G}{\overline P}<\infty$,
\item $\sdouble {gXg\inv }{G}{hP h\inv}<\infty$.
\end{theorem_count}
Let $G$ be an algebraic group and $X$ and $P$ be closed subgroups.
Denote by $G^\circ$, $X^\circ$ and $P^\circ$ the identity components
of $G$, $X$ and $P$ respectively.  If $\double XGP$ is finite then so
is $\double{X^\circ}{G^\circ}{P^\circ}$.
\end{lemma}

\begin{proof} 
These statements can all be proven in an elementary fashion.  The
final statement uses only the fact that $X^\circ$ and $P^\circ$ are
normal subgroups of finite index in $X$ and $P$ respectively.  
\end{proof}

This Lemma justifies assumptions implicit in the statement and in the
proof of Theorem~\ref{classification}.  The final statement of the
Lemma will be applied to recover finiteness results in $\SO(V)$ from
arguments made involving $\O(V)$.

\begin{remarks}
  We fix some notation for the rest of the paper.  We let $G$ be a
  classical algebraic group defined over a fixed algebraically closed
  field.  If $G$ is of type $B_n$, $C_n$ or $D_n$ we will, when
  convenient, assume that $G$ is one of $\SO_{2n+1}$, $\Sp_{2n}$ or
  $\SO_{2n}$ respectively.  We will often replace $G=A_n$ with
  $G=\GL_n$.  In all cases $G$ has rank $n$ (note that this means
  replacing $P_n$ in $A_n$ with $P_{n-1}$ in $\GL_n$).  We denote by
  $X$ a maximal rank reductive subgroup of $G$ and by $P$ a parabolic
  subgroup of $G$.
  
  We fix some terminology since usage varies in the literature and
  refer the reader to \cite{taylor} or \cite{grove} for further
  details.  To each possibility for $G$ we associate a natural module
  $V$ and a bilinear and quadratic form (we take these forms to be
  identically zero if $G=\GL_n$).  A subspace of $V$ is totally
  singular if each form is identically zero on the subspace.  If $P$
  is a parabolic subgroup then it equals the stabilizer of a (partial)
  flag of totally singular subspaces of $V$ and we identify $G/P$ as a
  collection of flags of totally singular subspaces.  Thus, a
  spherical subgroup has a finite number of orbits on the set of flags
  of totally singular subspaces.

Let $G$ equal $\SO_{2n}$.  Then two totally singular $n$-spaces are
conjugate under $G$ if and only if their intersection has odd
codimension in each space.
Thus, there exist two $G$-classes of totally singular $n$-spaces which
we identify as $G/P_{n-1}$ and $G/P_n$.  If $n$ is odd then $L_{n-1}$
is conjugate to $L_n$ under $G$, but this is not the case if $n$ is
even.  These facts are relevant to \Theorempart{classification}{(iv)}.

Let $V$ be the natural module for $G$.  We write $V= V_1 \perp V_2$ if
$V=V_1 \oplus V_2$ and each element of $V_1$ is orthogonal to each
element of $V_2$.  Given such a decomposition we define $\Cl(V_i)$ to
equal $\GL(V_i)$, $\O(V_i)$ or $\Sp(V_i)$ as $G$ equals, respectively,
$\GL(V)$, $\SO(V)$ or $\Sp(V)$.  We let $\Cl(V_i)^\circ$ be the
connected identity component of $\Cl(V_i)$.  If $\Cl(V_i)= \O(V_i)$
then $\Cl(V_i)^\circ = \SO(V_i)$ and otherwise $\Cl(V_i)^\circ =
\Cl(V_i)$.
\end{remarks}

\section{Finiteness}
\label{section_finiteness}
In this section we will prove those parts of
Theorem~\ref{classification} which assert finiteness.  For the
convenience of the reader we list in
Table~\ref{table_cases_finiteness} the specific result which covers
each case.
\begin{table*}
\caption{Finiteness Cases}%
\label{table_cases_finiteness}%
$$\begin{array}{cll}
\hline
G & \text{Cases} & \text{Proof}\\
\hline
\hline
A_n & \text{all cases} & \text{Corollary }\ref{finiteness_easy}\\
B_n & \text{the case with }P=P_1 & \text{Corollary }\ref{finiteness_easy}\\
    & \text{the case with }P=P_n & \text{Corollary }\ref{finiteness_of_AB_B_Pn}\\
C_n & \text{all cases with } P=P_1 & \text{Corollary }\ref{finiteness_easy}\\
    & \text{all cases with } P=P_n & \text{Corollary }\ref{finiteness_max_spaces}\\
D_n 
 & G=D_4,\text{ all cases} & \text{Corollary }\ref{cor_last_case_in_D4}\\
 & (X,P)=(A_{n_1}D_{n_2}T_1,\ P_1) & \text{Corollary }\ref{finiteness_easy}\\
 & (X,P)=(A_{n_1}A_{n_2}T_2,\ P_1) & \text{Lemma }\ref{finiteness_two_Ans_Dn_P1}\\
 & (X,P)\in \{D_{n_1}D_{n_2}D_{n_3},\ A_{n_1}D_{n_2}T_1\} \times 
            \{P_{n-1},\ P_n\} & \text{Corollary }\ref{finiteness_max_spaces}\\
\hline
\end{array}$$
\end{table*}







The next Lemma is not used immediately,
but we place it here to preserve the line of argument later.

\begin{lemma}
\label{lemma_GLn_on_1_spaces_in_SOn}
Let $G=\SO_{2n}$ and let $X=\GL_n$.  Then $X$ acts transitively upon
the set of definite $1$-spaces in $V$ and has a finite number of
orbits upon the set of all $1$-spaces in $V$.
\end{lemma}

\begin{proof}
  Let $N_1$ be the stabilizer of a definite $1$-space in $G$.  Then
  the first claim is equivalent to having $G=XN_1$.  In this form the
  first claim is proven in \cite{LSS}.  The second claim follows from
  the facts that $\GL_n$ is spherical in $G$ and that every $1$-space
  is either singular or definite.
\end{proof}

\begin{remark}
\label{notation_for_decomposition}
We fix some notation for the next few Lemmas.  Let $G$ be one of
$\GL(V)$, $\O(V)$, $\SO(V)$, or $\Sp(V)$ and fix a decomposition
$V=V_1\perpplus V_2$ such that $X= X_1X_2$ for subgroups $X_i \le
\Cl(V_i)$.  Let $P$ be a maximal parabolic subgroup of $G$ and
identify $G/P$ as a collection of totally singular subspaces in $V$.
For $i\in\{1,\ 2\}$ let $\pi^{}_i:V \to V_i$ be the natural projection.
\end{remark}

\begin{lemma} 
\label{lemma_projections_and_forms}
Let the notation be as in Remark~\ref{notation_for_decomposition}.  In
addition let $\beta$ and $\phi$ be, respectively, the bilinear and
quadratic forms associated with $G$.
\begin{theorem_count}
\item Let $u,\,v,\,x,\,y \in V$ such that $\beta(u,v)=\beta(x,y)=0$.  Then
$\beta(\pi^{}_1 u,\,\pi^{}_1 v)=\beta(\pi^{}_1 x,\, \pi^{}_1 y)$ if and only if
$\beta(\pi^{}_2 u,\,\pi^{}_2 v)=\beta(\pi^{}_2 x,\,\pi^{}_2 y)$.

\item Let $u$ and $x$ be two singular vectors with $\phi(\pi^{}_1 u) =
\phi(\pi^{}_1 x)$.  Then $\phi(\pi^{}_2 u)= \phi(\pi^{}_2 x)$.
\end{theorem_count}
\end{lemma}

\begin{proof}
  Decompose $u$ as $u=\pi^{}_1 u + \pi^{}_2 u$, decompose $v,\,x,$ and
  $y$ similarly, and use the fact that $V_1$ and $V_2$ are orthogonal
  to each other.
%
\end{proof}

\begin{lemma}
\label{lemma_reduction_to_X_1}
Let the notation be as in Remark~\ref{notation_for_decomposition} with
the additional assumptions that $G=\Cl(V)$ and $X_2 = \Cl(V_2)$.  Two
totally singular subspaces of the same dimension are conjugate under
$X$ if and only if their projections to $V_1$ and intersections with
$V_1$ are simultaneously conjugate under $X_1$.
\end{lemma}

The proof of this statement uses Witt's Theorem applied to $X_2$, so,
if $G$ is an orthogonal group, one cannot replace $\Cl(V_2) = \O(V_2)$
with $\Cl(V_2)^\circ = \SO(V_2)$.  We will use
Lemma~\ref{lemma_prelimnary_reductions} to translate finiteness
results to $\SO_n$.  We note that finiteness results do not always
translate between $\O_n$ and $\SO_n$ in an obvious fashion.  For
example, the collection $\double {L_{2,3}} {\SO_8} {P_4}$ is finite
whereas $\double {L_{2,3}} {\O_8} {P_4}$ is infinite.

\begin{proof} 
  It is easy to see that if two subspaces are conjugate under $X$,
  then their projections to $V_1$ and intersections with $V_1$ are
  simultaneously conjugate under $X_1$.
  
  Conversely, let $W$ and $W'$ be totally singular subspaces of the
  same dimension such that $x_1 ( W\cap V_1,\, \pi^{}_1 W) = ( W' \cap
  V_1,\,\pi^{}_1 W')$ for some $x_1\in X_1$.  Replacing $W$ with $x_1
  W$ we may assume that $ (W\cap V_1,\,\pi^{}_1 W) = (W'\cap
  V_1,\,\pi^{}_1 W')$.  Note that $\dim W\cap V_2=\dim W'\cap V_2$.
  Define the following dimensions:
\begin{eqnarray*} 
a & = & \dim W\cap V_1 = \dim W'\cap V_1,\\
c & = & \dim W\cap V_2 = \dim W\cap V_2,\\
b & = & \dim W-a-b = \dim W'-a-b.
\end{eqnarray*}
We will pick bases for $W$ and $W'$ as follows:
$$\begin{array}{rccc}
W: & w_1,\dots,w_{a}, & w_{a+1}, \dots, w_{a {+} b}, & w_{a {+} b{+}1},\dots,w_{a{+}b{+}c},\\
W': & w_1,\dots,w_{a}, & w_{a+1}',\dots,w_{a{+} b}', & w_{a{+}b+1}',\dots,w_{a{+}b{+}c}'.\\
\end{array}$$
We start by picking a basis $w_1,\dots,w_{a}$ of $W\cap V_1 = W'\cap
V_1$.  Extend this with elements $v_{a+1},\dots, v_{a+b}$ to a basis
of $\pi^{}_1 W = \pi^{}_1 W'$.  For each $i\in \{a+1,\dots,a+b\}$ pick
$w_i\in W$ and $w_i'\in W'$ such that $\pi^{}_1 w_i = \pi^{}_1
w_i'=v_i$.  Let $w_{a+b+1},\dots,w_{a+b+c}$ and $w_{a+b+1}', \dots,
w_{a+b+c}'$ be bases for $W\cap V_2$ and $W'\cap V_2$ respectively.
Then each of $\{w_{a+1},\dots,w_{a+b+c}\}$ and
$\{w_{a+1}',\dots,w_{a+b+c}'\}$ is a linearly independent set, whence
each of $\{\pi^{}_2 w_{a+1},\dots, \pi^{}_2 w_{a+b+c}\}$ and
$\{\pi^{}_2 w'_{a+1},\dots, \pi^{}_2 w_{a+b+c}'\}$ is a linearly
independent set.

Let $\tilde x_2$ be the linear map from the subspace $\genby{\pi^{}_2
  w_{a+1}, \dots, \pi_2^{} w_{a+b+c}}$ to $\genby{\pi^{}_2 w'_{a+1},
  \dots, \pi_2^{} w_{a+b+c}'}$ which takes each $\pi^{}_2 w_i$ to $
\pi^{}_2 w_i'$.

If $G=\GL(V)$ then one may extend $\tilde x_2$ to an element $x_2\in
\GL(V_2) = X_2$.  Note that $x_2 w_i = w_i'$ for each $i\in \{a+1, \dots,
a+b+c\}$.  This finishes the proof for the case $G=\GL(V)$.

If $G\in \{\O(V),\ \Sp(V)\}$ we show that $\tilde x_2$ is an isometry
from the subspace $\genby{\pi^{}_2 w_{a+1}, \dots, \pi_2^{}
  w_{a+b+c}}$ to $\genby{\pi^{}_2 w'_{a+1}, \dots, \pi_2^{}
  w_{a+b+c}'}$.  Once this is done, Witt's Theorem implies that we may
again extend $\tilde x_2$ to $x_2\in \Cl(V_2) = X_2$ and we will be
finished.  If $p\ne 2$ or if $G$ is symplectic, then
\Lemmapart{lemma_projections_and_forms}{(i)} shows that
$\beta(\pi_2^{} w_i, \pi_2^{} w_j) = \beta(\pi_2^{} w_i', \pi_2^{}
w_j')$ for all $a+1 \le i,j\le a+b+c$, whence $\tilde x_2$ is an
isometry.  If $p=2$ and $G$ is orthogonal, then
\Lemmapart{lemma_projections_and_forms}{(ii)} shows that
$\phi(\pi^{}_2w_i)=\phi(\pi^{}_2 w_i')$ for $a+1 \le i \le a+b+c$,
whence $\tilde x_2$ is an isometry.
\end{proof}

\begin{cor}
\label{main_cor_finiteness_classical}
Let the notation be as in Remark~{\rm\ref{notation_for_decomposition}}
with the additional assumption that $X_2$ equals $\Cl(V_2)$ or
$\Cl(V_2)^\circ$. Then $\double XGP$ is finite in the following cases:
\begin{theorem_count}
  
\item $X_1$ has a finite number of orbits upon the set of $1$-spaces
  in $V_1$ and $P=P_1$,
  
\item $G=\GL(V)$, $X_1$ has a finite number of orbits upon the set of
  all flags in $V_1$ and $P=P_i$ for some $i$,
  
\item $G = \GL(V) = \GL_n$, $X_1$ has a finite number of orbits upon
  the set of subspaces of $V_1$ with codimension $1$ and $P=P_{n-1}$.

\end{theorem_count}
\end{cor}

\begin{proof}
  By Lemma~\ref{lemma_prelimnary_reductions} it suffices to prove
  Corollary~\ref{main_cor_finiteness_classical} with the assumption
  that $G=\Cl(V)$ and $X=\Cl(V_2)$.  By
  Lemma~\ref{lemma_reduction_to_X_1} it suffices to show that $X_1$
  has a finite number of orbits upon the set $\{(W\cap V_1,\,\pi_1 W)
  \st W\in G/P\}$.  Note that $W\cap V_1 \le \pi_1 W$ is a flag.  It
  is easy to verify in each case that $X_1$ has a finite number of
  orbits on the set.
\end{proof}

\begin{cor}
\label{finiteness_easy}
The double coset collection $\double XGP$ is finite in the following
cases:
\begin{theorem_count}
  
\item $G=\GL_{n}$, {\rm (a)} $X$ has no additional restrictions (other
  than our standing assumptions) and $P\in \{P_1,P_{n-1}\}$ or {\rm
    (b)} $X=\GL_{n_1}\GL_{n_2}\GL _{n_3}$ and $P=P_i$ for some $i$,

\item $G=\SO_{2n+1}$, $X=\GL_{n_1}\SO_{2n_2+1}$ and $P = P_1$,
  
\item $G=\Sp_{2n}$, $X\in\{\GL_{n_1}\Sp_{2n_2}\cdots\Sp_{2n_r},\ 
  \Sp_{2n_1}\cdots\Sp_{2n_r}\}$ and $P=P_1$,

\item $G= \SO_{2n}$, $X= \GL_{n_1}\SO_{2n_2}$ and $P = P_1$.

\end{theorem_count}
\end{cor}

\begin{proof}
In each case let $V$ be the natural module for $G$ and fix a
decomposition $V = V_1 \perp V_2$ such that $X=X_1X_2$ with $X_1 \le
\Cl(V_1)$ and $X_2 = \Cl(V_2)^\circ$.  

For case (i)(b) note that $X_1= \GL_{n_1}\GL_{n_2}$ is a spherical
subgroup of $\GL_{n_1+n_2} = \Cl(V_1)$ and apply
\Corpart{main_cor_finiteness_classical}{(ii)}.  For case (i)(a) apply
\Corpart{main_cor_finiteness_classical}{(i)} if $P=P_1$, apply 
\Corpart{main_cor_finiteness_classical}{(iii)} if $P=P_{n-1}$, and, in
both cases, induct on $r$.

Let $G=\Sp_{2n}$.  By \Corpart{main_cor_finiteness_classical}{(i)} it
suffices to show that $X_1$ has finitely many orbits on totally
singular $1$-spaces in $V$ (note that all $1$-spaces are totally
singular in this case).  This is immediate if $X_1$ equals $\GL_{n_1}$
or $\Sp_{2n_1}\Sp_{2n_2}$ since these subgroups are spherical in
$\Cl(V_1)^\circ$.  The general case follows by induction on $r$.

If $G$ is orthogonal, then $X_1=\GL_{n_1}$ and the result follows by
combining Lemma~\ref{lemma_GLn_on_1_spaces_in_SOn} and
\Corpart{main_cor_finiteness_classical}{(i)}.
\end{proof}

\begin{lemma} 
\label{symplectic_lemma:lasso_maximal}
Let the notation be as in Remark~\ref{notation_for_decomposition} with
the additional assumption that $V$ is a symplectic or orthogonal
space.  Let $W$ be a maximal totally singular subspace of $V$ and let
$(\pi_i^{}W)^\perp$ be the perpendicular space taken within $V_i$.  If
$\dim V$ is even then $\dim (\pi_i^{} W)^\perp
/(W\cap V_i)$ equals $0$.  If $\dim V$ is odd then $\dim (\pi_i^{} W)^\perp
/(W\cap V_i)$ equals $0$ or $1$.
\end{lemma}

\begin{proof}
Since $V_1$ is orthogonal to $V_2$ it is easy to show that $W\cap V_i
\le (\pi_i W)^\perp$.
We have $\dim W$ equals $n$ and $\dim V$ equals $2n$ or $2n+1$.  Set
$a_i = \dim W\cap V_i$.  For $i$ equal to $1$ and $2$ the inequality
$\dim (\pi_i W)^\perp \ge \dim W\cap V_i$ becomes, respectively, $\dim
V_1-(n-a_2)\ge a_1$ and $\dim V_2-(n-a_1)\ge a_2$.  If $\dim V$ is
even, then the sum of these last two inequalities is an equality; if
$\dim V$ is odd, then the sum of the left sides is 1 greater than the sum
of the right sides.
\end{proof}

\begin{cor}   
\label{corollary:action_on_max}
Let the notation be as in Remark~\ref{notation_for_decomposition} with
the additional assumptions that $G$ is not $\GL(V)$, that $X_2$ equals
$\Cl(V_2)^\circ$ or $\Cl(V_2)$ and that $P=P_n$.  Then $\double XGP$
is finite in the following cases:
\begin{theorem_count}
  
\item $\dim V$ is even and $X_1$ has a finite number of orbits on the
  set of totally singular subspaces of $V_1$,

\item $\dim V$ is odd and $X_1$ has a finite number of orbits on the
  set of pairs of subspaces $(W_1,\,W_2)$ such that $W_1 \le W_2 \le
  V_1$, $W_1$ is totally singular, $W_1\le W_2^\perp$ and $\dim
  (W_2/W_1)\le 1$.
\end{theorem_count}
\end{cor}

\begin{proof} 
  By Lemma~\ref{lemma_prelimnary_reductions} we may assume that
  $G=\Cl(V)$ and $X_2=\Cl(V)$.  By Lemma~\ref{lemma_reduction_to_X_1}
  it suffices to show that $X_1$ has a finite number of orbits upon
  the set $\{(W\cap V_1, \pi_1^{} W) \st W\in G/P\}$.  Given subspaces
  $W,\,W'\le V$ and $x_1 \in X_1$ we have $x_1 \pi_1^{}W = \pi_1^{}
  W'$ if and only if $x_1 (\pi_1^{} W)^\perp = (\pi_1^{}W')^\perp$,
  where we take the perpendicular space within $V_1$.  Thus, it
  suffices to show that $X_1$ has a finite number of orbits on the set
  $\{(W\cap V_1,\ (\pi_1^{} W)^\perp) \st W \in G/P\}$.  By
  Lemma~\ref{symplectic_lemma:lasso_maximal} we see that $\{(W\cap
  V_1,\ (\pi_1^{} W)^\perp) \st W\in G/P\}$ is a subset (or may be
  identified with a subset) of one of the sets given in the statement
  of Corollary~\ref{corollary:action_on_max}.
%
\end{proof}

\begin{cor}
\label{finiteness_max_spaces}
The double coset collection $\double XGP$ is finite in
the following cases:
\begin{enumerate}
\compress\romanate

\item $G=\SO_{2n}$, $X\in \{\SO_{2n_1}\SO_{2n_2}\SO_{2n_3}$,
  $\GL_{n_1}\SO_{2n_2}\}$ and $P\in \{P_{n-1},\ P_n\}$.

\item $G=\Sp_{2n}$, $X\in \{\Sp_{2n_1}\Sp_{2n_2}\Sp_{2n_3}$,  
$\GL_{n_1}\Sp_{2n_2}$,  $\Sp_{2n_1}\Sp_{2n_2}T_1\}$ and $P=P_n$. 

\end{enumerate}
\end{cor}

\begin{proof}
  Let $V$ be the natural module of $G$ and fix a decomposition $V=V_1
  \perp V_2$ so that $X=X_1X_2$ with $X_1 \le \Cl^\circ (V_1)$ and
  $X_2=\Cl(V_2)^\circ$.  Then $X_1$ is a spherical subgroup of
  $\Cl(V_1)^\circ$
whence the conclusion follows from
\Corpart{corollary:action_on_max}{(i)}.
\end{proof}

\begin{lemma}
\label{descent_in_Dn} 
Let $G = \SO_{2n}$, let $X= \GL_n\le \SO_{2n}$ and let $V$ be the
natural module for $G$.  Let $\genby v$ be a definite $1$-space, let
$X_{\genby v}$ be the stabilizer in $X$ of $\genby v$ and let $\tilde
X_{\genby v}$ be the connected component of the group induced by
$X_{\genby v}$ in $\SO({\genby v}^\perp)$.  Then $\tilde X_{\genby v}$
is a spherical subgroup of $\SO({\genby v}^\perp)$.
\end{lemma}

\begin{proof}
  By Lemma~\ref{lemma_GLn_on_1_spaces_in_SOn} we may calculate
  $X_{\genby{v}}$ where $v$ is any definite vector.
  Let $V = E\oplus F$ such that $E$ and $F$ are totally singular and
  $X$ is the stabilizer in $G$ of this decomposition.  Let $\{e_1,
  \dots, e_n\}\subset E$ and $\{f_1, \dots, f_n\} \subset F$ be dual
  bases.  Let $v=e_1+f_1$ and let $\GL_{n-1}$ denote the subgroup of
  $X$ which stabilizes the subspaces $\genby{e_2, \dots, e_n}$ and
  $\genby{f_2, \dots, f_n}$ and acts trivially upon $e_1$ and $f_1$.
  Then (an isomorphic image of) $\GL_{n-1}$ is a subgroup of $\tilde
  X_{\genby v}$ which proves the claim since $\SO(\genby v^\perp) =
  \SO_{2(n-1)+1}$.

%
\end{proof}

%
%

\begin{cor} 
\label{finiteness_of_AB_B_Pn}
Let $G=\SO_{2n+1}$, $X=\GL_{n_1}\SO_{2n_2+1}$, and $P=P_n$.  Then
$\double XGP$ is finite.
\end{cor}

\begin{proof}
  Let $V$ be the natural module of $G$ and fix a decomposition $V=V_1
  \perp V_2$ so that $X=X_1X_2$ with $X_1 = \GL_{n_1} \le \Cl^\circ
  (V_1)$ and $X_2= \SO_{2n_2+1} = \Cl(V_2)^\circ$.  By
  \Corpart{corollary:action_on_max}{(ii)} it suffices to show that
  $X_1$ has a finite number of orbits on the set of pairs of subspaces
  $(W_1,\,W_2)$ such that $W_1 \le W_2 \le V_1$, $W_1$ is totally
  singular, $W_1 \le W_2^\perp$ and $\dim (W_2/W_1)\le 1$.  We may
  partition this set into two subsets according as $W_2$ is, or is
  not, totally singular.  Since $X_1$ is spherical in $\SO(V_1)$ we
  see that $X_1$ has finitely many orbits upon the subset where $W_2$
  is totally singular.
  
  Every pair $(W_1,\,W_2)$ where $W_2$ is not totally singular can be
  rewritten as $(W_1,\, W_1{\perp} \genby{v})$ where $v$ a definite
  vector in $V_1$.  Thus it suffices to show that $X_1$ has a finite
  number of orbits upon pairs $(W_1,\genby{v})$ such that $v$ is a
  definite vector in $V_1$ and $W_1 \le \genby{v}^\perp$ is totally
  singular (where this perpendicular space is taken in $V_1$).  By
  Lemma~\ref{lemma_GLn_on_1_spaces_in_SOn}, $X_1$ acts transitively
  upon definite $1$-spaces, whence it suffices to fix $v$ and show
  that the stabilizer in $X_1$ of $\genby v$ has a finite number of
  orbits on totally singular subspaces in $\genby{v}^\perp$.  This
  follows from Lemma~\ref{descent_in_Dn}.
\end{proof}

\begin{lemma} 
\label{finiteness_two_Ans_Dn_P1}
Let $G=\SO_{2n}$, $X=\GL_{n_1}\GL_{n_2}$ and $P=P_1$.  
Then $\double XGP$ is finite.
\end{lemma}

\begin{proof}
  Let $V$ be the natural module of $G$ and fix a decomposition $V=V_1
  \perp V_2$ so that $X=X_1X_2$ with $X_i = \GL_{n_i} \le \Cl^\circ
  (V_i)$ for each $i$.
  
  Let $\pi_i^{}:V\to V_i$ be the natural projection.  We wish to show
  that $X$ has a finite number of orbits on the set $\{\genby{v} \st v
  \in V \text{ is singular}\}$.  By
  Lemma~\ref{lemma_GLn_on_1_spaces_in_SOn} each $X_i$ has finitely
  many orbits on $1$-spaces in $V_i$.  Thus it suffices to show, for
  each singular $1$-space $\genby{v}$ that $X$ has finitely many
  orbits on the set of singular $1$-spaces whose projections to $V_1$
  and $V_2$ are conjugate to $\pi_1^{}\genby{v}$ and
  $\pi_2^{}\genby{v}$ respectively.  Thus, it suffices to fix an
  arbitrary singular $1$-space $\genby{v}$, let $v_i=\pi_i^{}v$ and
  show that $X$ has a finite number of orbits on the set
\begin{equation}
\label{set_of_skewed_1_spaces}
\{\genby{v_1+\alpha v_2} \st \alpha\in \k,\ \alpha\ne0,\ v_1+\alpha
v_2 \text{ is singular}\},
\end{equation}
where $\k$ is the ground field.  Since $v_1 +\alpha v_2$ is singular,
we have $\phi(v_1) + \alpha ^2 \phi(v_2)=0$, where $\phi$ is the
quadratic form.  Thus, $v_1$ is definite if and only if $v_2$ is.  If
$v_1$ and $v_2$ are both definite then $\alpha^2 =
-\phi(v_1)/\phi(v_2)$, whence there are at most two singular
$1$-spaces of the form $\genby {v_1+\alpha v_2}$.

It suffices now to assume that $v_1$ and $v_2$ are singular and show
that $X_1$ has finitely many orbits on the set in
Equation~\ref{set_of_skewed_1_spaces}.  Let $V_2= E_2\oplus F_2$ where
$E_2$ and $F_2$ are totally singular and $X_2$ stabilizes $E_2$ and
$F_2$.  Then $v_2= e + f$ for some $e_2\in E_2$, $f_2\in F_2$.  Since
$v_2$ is singular this implies that $\beta(e,f)=0$, where $\beta$ is
the bilinear form.  Easy linear algebra shows that there exists $x_2
\in X_2$ with $x_2 e = \alpha e$ and $x_2 f = \alpha f$, whence $x_2
\genby{v_1+v_2} = \genby{v_1+\alpha v_2}$.
\end{proof}

Recall that $L_{i,2}$ denotes a Levi subgroup as described just before 
Theorem~\ref{classification}.  

\begin{cor}
  \label{cor_last_case_in_D4}
  Let $G = D_4$ and $(X,P) \in \{(L_{2,3},\ P_4),\ (L_{2,4},\ P_{3})\}$.
  Then $\double XGP$ is finite.
\end{cor}

\begin{proof}
  This follows most easily from applying the graph automorphism of
  order three to other cases which have been proven finite.  For
  instance $(X,P) = (L_{2,3}, P_4)$ follows from
  Corollary~\ref{finiteness_max_spaces} applied to $(X,P) = (L_{1,2},
  P_3)$ (with $L_{1,2} = D_1D_2D_2$) or from
  Lemma~\ref{finiteness_two_Ans_Dn_P1} applied to $(X,P) = (L_{2,4},
  P_1)$.
\end{proof}

The reader who wishes for an instructive, though somewhat painful,
exercise can prove this Corollary using geometric arguments about
subspaces of the natural module of $D_4$.


\section{Infiniteness}
\label{infinite}

We begin by stating a result which gives infiniteness in many cases.

\begin{theorem}[\cite{duckworth}*{Theorem 1.3}]
\label{sph_or_sph_levi}
If $\double XGP$ is finite then $X$ or $L$ is a spherical subgroup of
$G$.
\end{theorem}

\begin{remark}
  To finish the proof of infiniteness in Theorem~\ref{classification}
  it suffices, by Theorem~\ref{sph_or_sph_levi}, to consider only
  those $P$ such that $L$ is spherical.  Suppose that we have fixed such a
  $P$.  Then it suffices to prove infiniteness for those $X$ which are
  maximal subject to the condition that $\double XGP$ is claimed to be
  infinite in Theorem~\ref{classification}.  
  
  In Table~\ref{table_for_infiniteness_proof} we list those cases
  which need to be proven infinite, and indicate which result
  addresses each case.  Recall that we allow the notation $A_0$,
  $B_0$, $C_0$ and $D_1$ unless otherwise noted (but we do not allow
  $D_0$).
\end{remark}

\begin{table*}
\caption{Infiniteness cases}%
\label{table_for_infiniteness_proof}%
$$
\begin{array}{lll}
\hline 
G & (X,P) & \text{Proof}\\
\hline\hline
A_n 
 & \{A_{n_1}A_{n_2}A_{n_3}A_{n_4}T_3\} \times \{P_i \st 2\le i\le n-1\} 
      & \text{Lemma }\ref{lemma_An_GLnGLnGLnGLn_Pi}\\
B_n
 & \{B_{n_1}D_{n_2}D_{n_3},\ A_{n_1}A_{n_2}T_2\}\times \{P_1,\ P_n\} 
     & \text{Lemma }\ref{lemma_Bn_all_cases}\\
C_n
 & \{A_{n_1}A_{n_2}C_{n_3}T_2\} \times \{P_1,\ P_n\}
 & \text{Lemma }\ref{lemma_Cn_all_cases} \\
 & (C_{n_1}C_{n_2}C_{n_3}C_{n_4}\ (n_i\ge 1),\ P_n)
    & \text{Lemma }\ref{lemma_Cn_all_cases} \\
 & (A_{n_1}C_{n_2}C_{n_3}T_1\ (n_i\ge 1),\ P_n)  
    & \text{Lemma }\ref{lemma_Cn_all_cases} \\
D_n
 & (D_{n_1}D_{n_2}D_{n_3},\ P_1) 
   & \text{Lemma }\ref{lemma_Dn_all_cases_except_AnAn}   \\
 & \{D_{n_1}D_{n_2}D_{n_3}D_{n_4}\} \times \{P_{n-1},\ P_n\} 
   & \text{Lemma }\ref{lemma_Dn_all_cases_except_AnAn}  \\ 
 & \{ A_{n_1}D_{n_2}D_{n_3}\ (n_1\ge1)\} \times \{P_{n-1},\ P_n\} 
   & \text{Lemma }\ref{lemma_Dn_all_cases_except_AnAn}  \\[3\jot]
 & \{ A_{n_1}A_{n_2}T_2\ (n_i \ge 1)\}  \times \{P_{n-1},\ P_n\} &
   \text{Lemma }\ref{lemma_Dn_GLnGLn_Pn} \\
 & \text{with } (G,X,P) \text{ not as in \Theorempart{classification}{(iv)(a)}}
   &  \\
\hline 
\end{array}
$$
\end{table*}

\begin{remark}  For the remainder of this section we assume
  that $X$ and $L$ contain a common maximal torus, $T$.  For a closed
  subgroup $H$ which contains $T$ we write $\rtsys(H)$ for the root
  system of $H$ defined using $T$.  For a closed root subsystem $\phi$
  of $\rtsys(G)$ we set $\dim \phi = \dim H/Z$ where $H$ is a closed
  subgroup of $G$ such that $\rtsys(H) = \phi$ and $Z$ is the center
  of $H$ (Theorem~\ref{theorem_root_complements} also holds if we use
  $\dim \phi = \dim H$ instead).
\end{remark}

\begin{theorem}[\cite{duckworth}*{Theorem 1.1, Lemma 3.3}]
  \label{theorem_root_complements}
  For $i\in \{1,2\}$ let $L_i$ be a Levi subgroup containing $T$ and
  $\rtsys_i$ its root system.  Assume that $L_1$ and $L_2$ are
  conjugate.  If $\frac 12 \dim \rtsys_1 - \dim \rtsys_1\cap \rtsys(X)
  - \frac 12 \dim \rtsys_2 \cap \rtsys(L) > 0$ then $\double XGP$ is
  infinite.  In particular infiniteness holds in the following cases:
\begin{theorem_count}
\item $\rtsys_1$ and $\rtsys_2$ are of type $B_2$, $\rtsys_1\cap
      \rtsys(X) = \emptyset$ and $\rtsys_2 \cap \rtsys(L)$ is of type
      $A_1$.

\item $\rtsys_1$ and $\rtsys_2$ are of type $A_3$ or $D_3$, $\rtsys_1
      \cap \rtsys(X) = \emptyset$ and $\rtsys_2$ is of type $A_1A_1$
      or $D_2$.  
\end{theorem_count}
\end{theorem}

\begin{remarks} 
  We offer comments which help simplify the proofs of Lemmas
  \ref{lemma_An_GLnGLnGLnGLn_Pi}, \ref{lemma_Bn_all_cases},
  \ref{lemma_Cn_all_cases} and \ref{lemma_Dn_all_cases_except_AnAn}.
  
  \rem All Levi subgroups of type $B_2$ are conjugate and, unless
  $G=D_n$, all Levi subgroups of type $A_3$ are conjugate.  Thus, to
  apply Theorem~\ref{theorem_root_complements} one often only has to
  verify that $\rtsys_1$, $\rtsys_1\cap \rtsys(X)$, $\rtsys_2$ and
  $\rtsys_2\cap \rtsys(L)$ are of the required type.  We will
  construct each $\rtsys_i$ by giving a base $\alpha$, $\beta$, \dots\ 
  and setting $\rtsys_i$ equal to all the $\Z$-linear combinations of
  $\alpha$, $\beta$, \dots\ which are in $\rtsys(G)$.  \rem Let
  $\rtbase(G)$ and $\tilde \rtbase(G)$ be the Dynkin diagram and
  extended Dynkin diagram of $G$ respectively.  Label the nodes of
  $\rtbase(G)$ using $\alpha_1$, \dots, $\alpha_n$ as in
  \cite{bourbaki}.  In each of the following Lemmas we will assume
  that $\rtbase(X)$ has been produced from $\rtbase(G)$ by the
  Borel-de Siebenthal algorithm \cite{borel_siebenthal}.  Let
  $X=X_{n_1}X_{n_2}\cdots$ with each $X_{n_i}$ equal to $D_1$, $T_1$
  or a simple factor of $X$, as listed in
  Table~\ref{table_for_infiniteness_proof}, with $n_i$ the rank of the
  factor.  We take $\rtbase(X_{n_1})$ equals to $\emptyset$ if
  $X_{n_1}$ equals $T_1$ or $D_1$ and otherwise $\rtbase(X_{n_1}$
  equals the first $n_1$ nodes of $\rtbase(G_1)$ or $\tilde
  \rtbase(G_1)$ as appropriate.  We then repeat this procedure,
  starting with $\rtbase(X_{n_2})$ and the last $n_2+n_3 + \dots$\ 
  nodes of $\rtbase(G)$.  This procedure determines
  $\rtbase(G)-\rtbase(X)$ which, in turn, provides an easy description
  of $\rtsys(X)$.  For example, suppose that $G=D_n$ and $X =
  D_{n_1}D_{n_2}$.  Then $\rtbase(G)-\rtbase(X) = \{\alpha_{n_1}\}$
  and $\rtsys(X)$ equals all the roots in $\rtsys(G)$ which have
  $\alpha_{n_1}$-coefficient equal to $0$ or $\pm2$.  
  \rem Given $\alpha$, $\beta \in \rtbase(G)$ the path connecting
  $\alpha$ to $\beta$ is the shortest such path and includes $\alpha$
  and $\beta$.  The sum over this path means the sum of each element
  of $\rtbase(G)$ which is contained in the path.  It is easy to check
  that such a sum is itself a root.
\end{remarks}

\begin{lemma_case}
  \label{lemma_An_GLnGLnGLnGLn_Pi}
  Let $G = A_n$, $X = A_{n_1}A_{n_2}A_{n_3}A_{n_4}T_3$ and $P\in \{P_i
  \st 2\le i \le n-1\}$.  Then there exist $\rtsys_1$ and $ \rtsys_2$
  of type $A_3$ as in \Theorempart{theorem_root_complements}{(ii)}.
\end{lemma_case}

\begin{proof}
  We have $\rtbase(G) - \rtbase(X) = \{\alpha_{n_1{+}1}$,
  $\alpha_{n_1{+}n_2{+}2}$, $\alpha_{n_1{+}n_2{+}n_3{+}3}\}$.  Let
  $\rtsys_1$ have root base given by $\alpha$ equal to
  $\alpha_{n_1+1}$, $\beta$ equal to the sum over the path connecting
  $\alpha_{n_1+2}$ to $\alpha_{n_1+n_2+2}$, and $\gamma$ equal to the
  sum over the path connecting $\alpha_{n_1+n_2+3}$ to
  $\alpha_{n_1+n_2+n_3+3}$.
  
  For $L_i$ let $\rtsys_2$ have root base given by $\alpha =
  \alpha_{i-1}$, $\beta = \alpha_i$, and $\gamma = \alpha_{i+1}$.
\end{proof}


\begin{lemma_case}
\label{lemma_Bn_all_cases}
Let $G=B_n$ and $(X,P) \in \{B_{n_1}D_{n_2}D_{n_3},\ 
A_{n_1}A_{n_2}T_2\} \times \{P_1, P_n\}$.  Then there exist $\rtsys_1$
and $\rtsys_2$ of type $ B_2$ as in
\Theorempart{theorem_root_complements}{(i)}.
\end{lemma_case}

\begin{proof}
  One may proceed as in the proofs of Lemmas
  \ref{lemma_An_GLnGLnGLnGLn_Pi} and \ref{lemma_Dn_GLnGLn_Pn}, or as
  in \cite{duckworth}*{Corollary 7.2 (ii)}.
\end{proof}

\begin{lemma_case}
\label{lemma_Cn_all_cases}
Let $G=C_n$.  If $(X,P) \in \{A_{n_1}A_{n_2}C_{n_3}T_2\}
\times\{P_1,P_n\}$, then there exists $\rtsys_1 = \rtsys_2$ of type
$B_2$ as in \Theorempart{theorem_root_complements}{(i)}.  If $(X,P)
\in \{C_{n_1}C_{n_2}C_{n_3}C_{n_4}\ \linebreak[1](n_i\ge 1)$,
$A_{n_1}C_{n_2}C_{n_3}T_1\ (n_i\ge 1)\}\times\{P_n\}$, then there
exist $\rtsys_1$ and $\rtsys_2$ of type $A_3$ as in
\Theorempart{theorem_root_complements}{(ii)}.
\end{lemma_case}

\begin{proof}
  One may proceed as in the proofs of Lemmas
  \ref{lemma_An_GLnGLnGLnGLn_Pi} and \ref{lemma_Dn_GLnGLn_Pn}.  
\end{proof}

\begin{lemma_case}
\label{lemma_Dn_all_cases_except_AnAn}
Let $G$ equal $ D_n$.  If $(X,P)$ equals $(D_{n_1}D_{n_2}D_{n_3},\ 
P_1)$, or is in $\{A_{n_1}D_{n_2}D_{n_3}T_1\ (n_1\ge 1),\ 
D_{n_1}D_{n_2}D_{n_3}D_{n_4}\} \times \{P_{n-1},\ P_n\}$, then there
exists $\rtsys_1 = \rtsys_2$ of type $A_3$ or $D_3$ as in
\Theorempart{theorem_root_complements}{(ii)}.
\end{lemma_case}

\begin{proof}
  One may proceed as in the proofs of Lemmas
  \ref{lemma_An_GLnGLnGLnGLn_Pi} and \ref{lemma_Dn_GLnGLn_Pn}.%
\end{proof}

\begin{lemma_case}
  \label{lemma_Dn_GLnGLn_Pn}
  Let $G=D_n$, $X = A_{n_1}A_{n_2}T_2\ (n_i\ge1)$, $P\in\{P_{n-1},\ 
  P_n\}$.  If $n=4$ we assume that $(X,P)$ is not equal to either
  $(L_{2,3},\ P_4)$ or $(L_{2,4},\ P_3)$.  Then there exists
  $\rtsys_1$ and $\rtsys_2$ of type $A_3$ or $D_3$ as in
  \Theorempart{theorem_root_complements}{(ii)}.
\end{lemma_case}

\begin{proof}
  By our convention with subsystems of type $A_{n_i}$ in $D_n$, we
  have that $X\in \{L_{i,n{-}1},\ L_{i,n}\}$ where $i=n_1+1$ satisfies
  $2\le i \le n-2$.
  
  Let $(X,P) = (L_{i,n},\ P_n)$.  Let $\rtsys_1 = \rtsys_2$ have root
  base given by $\alpha$ equal to the sum over the path connecting
  $\alpha_1$ to $\alpha_{n-1}$, $\beta$ equal to $\alpha_n$, and
  $\gamma$ equal to the sum over the path connecting $\alpha_2$ to
  $\alpha_{n-2}$.

  By symmetry the conclusion holds also when $(X,P) = (L_{i,n{-}1},\
  P_{n{-}1})$.  This leaves the cases $(X,P) \in \{(L_{i,n{-}1},\ P_n),\ 
  (L_{i,n},\ P_{n{-}1})\}$.  If $n$ is odd then $L_{i,n{-}1}$ is conjugate
  to $L_{i,n}$, whence the conclusion holds.
  It remains to prove the existence of $\rtsys_1$ and $\rtsys_2$ when
  $n \ge 6$ is even.  If necessary we replace $X$ by a conjugate to
  assume that $i \le n-3$.
  Let $(X,P) = (L_{i,n{-}1}, P_n)$ with $2 \le i \le n-3$.  Let
  $\rtsys_1= \rtsys_2$ have root base given by $\alpha$ equal to the
  sum over the path connecting $\alpha_i$ to $\alpha_{n-3}$, $\beta =
  \alpha_{n-2} + \alpha_{n-1} + \alpha_n$, and $\gamma$ equal to the
  sum over the path connecting $\alpha_{i-1}$ to $\alpha_{n-2}$.  By
  symmetry the conclusion also holds when $(X,P) = (L_{i,n},
  P_{n{-}1})$.  
\end{proof}

\end{document}